\documentclass[12pt]{amsart}
\usepackage{bbm}
\usepackage{amscd, amssymb, dsfont, upgreek, mathrsfs}
\input{xy}
\xyoption{all}

\newtheorem{Thm}{Theorem}[section]

\newtheorem{Prop}[Thm]{Proposition}

\newtheorem{Def/Thm}[Thm]{Definition/Theorem}
\newtheorem{Cor}[Thm]{Corollary}

\theoremstyle{remark}

\numberwithin{equation}{section}
\newcommand{\ti }{\times}
\newcommand{\ot }{\otimes}
\newcommand{\ra }{\rightarrow}

\newcommand{\cO}{{\mathcal{O}}}

\newcommand{\cC}{{\mathcal{C}}}

\newcommand{\G}{{\bf G}}

\newcommand{\PP }{{\mathbb P}}

\newcommand{\QQ }{{\mathbb Q}}
\newcommand{\CC }{{\mathbb C}}
\newcommand{\ZZ }{{\mathbb Z}}

\newcommand{\ke }{{\varepsilon }}

\newcommand{\kl }{{\lambda}}

\newcommand{\vir}{\mathrm{vir}}

\newcommand{\T}{{\bf T}}

\newcommand{\re}{\mathrm{e}}

\newcommand{\lan}{\langle}
\newcommand{\ran}{\rangle}
\newcommand{\lla}{\langle\!\langle}
\newcommand{\rra}{\rangle\!\rangle}

\newcommand{\WtG}{W/\!\!/_{\!\theta}\G}

\newcommand{\bx}{\mathbf{x}}

\newcommand{\ret}{\re ^{\T}}

\begin{document}

\title[Mirror Theorem for Elliptic Quasimap Invariants]{Mirror Theorem for Elliptic Quasimap Invariants of Local Calabi-Yau Varieties}

\begin{abstract} 
The elliptic quasi-map potential function is explicitly calculated for Calabi-Yau complete intersections in projective spaces in \cite{KL}. We extend this result to local Calabi-Yau varieties. Using this as well as the wall crossing formula in \cite{CKw}, we can calculate the elliptic Gromov-Witten potential function.

\end{abstract}

\author{Hyenho Lho}
\noindent\address{
Korea Institute for Advanced Study,
85 Hoegi-ro, Dongdaemun-gu, Seoul, 02455, Korea}
\email{hyenho@kias.re.kr}

\author{Jeongseok Oh}
\noindent\address{
Korea Advanced Institute of Science and Technology,
291 Daehak-ro, Yuseong-gu, Daejeon 34141, Korea}
\email{ojs0603@kaist.ac.kr}


\maketitle


\section{Introduction}

For a nonsingular variety $X$ which has a GIT representation $\WtG$, we can define the moduli spaces of $\ke$-stable quasimaps with genus $g$, $k$-markings to $X$ with degree $\beta$, denoted by $Q^{\ke}_{g, k}(X, \beta)$, for any $g$, $k$ with $2g-2+k \geq 0$, $\beta \in \text{Hom}_\mathbb{Z}(\text{Pic}^{\textbf{G}}(W),\mathbb{Z})$ unless $2g-2+k=0$ and $\beta=0$. For each $Q^{\ke}_{g, k}(X, \beta)$, we can define the canonical virtual fundamental class 
$$[Q^{\ke}_{g, k}(X, \beta)]^{\vir} \in A_* (Q^{\ke}_{g, k}(X, \beta))\otimes_{\mathbb{Z}} \mathbb{Q}$$
of degree
$$c^{\textbf{G}}_1(W).\beta + (\text{dim}_\mathbb{C} X-3)(1-g) +k.$$
See \cite{CKM} for details.

Especially, for a Calabi-Yau variety $X$, since $c^{\textbf{G}}_1(W)=0$, every $[Q^{\ke}_{1, 0}(X, \beta)]^{\vir}$ for any $\beta \neq 0$ has degree $0$. So, we can define a generating function
 \begin{align*}  \lan\; \ran ^{\ke }_{1,0} := \sum_{\beta \neq 0} q^\beta \deg [Q^{\ke}_{1,0}(X, \beta)]^{\vir} \end{align*}
for each $\ke$.
In particular, when $\ke$ is small enough, i.e. $\ke =0+$, it is called the elliptic quasi-map potential function of $X$.

Throughout this paper, let $X$ be a total space of vector bundle 
$$\mathcal{O}_{\PP^{n-1}}(-\l'_1)|_{X'}\oplus \mathcal{O}_{\PP^{n-1}}(-\l'_2)|_{X'}\oplus\dots \oplus \mathcal{O}_{\PP^{n-1}}(-\l'_m)|_{X'}$$ 
over $X'$, where $X'$ is a complete intersection in $\PP ^{n-1}$ defined by deg $l_i$ polynomials for $i=1,2,\dots,r$ and $l_a$, $l'_b>0$ for all $a$, $b$. 
We assume Calabi-Yau condition 
$$\sum_a l_a + \sum_b l'_b = n. $$
Note that $X$ has a natural GIT representation and is a Calabi-Yau variety.
In this paper, we will give an explicit formula of elliptic quasi-map potential function for this $X$.
In \cite{KL}, they already computed elliptic quasi-map potential function for $m=0$ case. Basically, we will follow their idea to prove the main theorem which we will introduce soon except for computational part. That is the following:
By Quantum Lefschetz hyperplane section theorem \cite{Kim}, quasi-map invariants of X can be represented as twisted quasi-map invariants of $\PP^{n-1}$ which also has a natural GIT representation. And we will apply the torus localization theorem for the later since $\PP^{n-1}$ has a natural torus action.

In order to state the main theorem, we need to give some preparations first.
In \cite{Gequiv}, Givental introduced equivariant I-function for X which is $H^*_{ \T }(\PP ^{n-1})\ot \QQ(\lambda, \zeta)$-valued formal function in formal variables $q,z,t_H$ :
$$ I^\zeta_{\T}(t,  q):= e^{t_HH/z}\sum _{d=0}^{\infty} q^d e^{t_Hd}
 \frac{\prod _{j=1}^r  \prod _{k=1}^{l_jd} (l_jH + kz)\prod _{j=1}^m  \prod _{k=0}^{l'_jd-1} (-l'_jH- kz+\zeta)}{\prod _{k=1}^{d} \prod _{j=1}^n (H-\lambda _j +kz)} , $$
where $\T =(\CC^*)^n$ is the torus group acting on $\PP^{n-1}$; $\lambda _1, ...,\lambda _n$ are the $\T$-equivariant parameter; $\zeta$ is the $\CC^*$-equivariant parameter for $\CC^*$-action acting diagonally on the fiber of $X$ over $X'$. $\QQ(\lambda, \zeta)$ denotes the quotient field of the polynomial ring in $\lambda _1, ...,\lambda _n, \zeta$; $H$ is the hyperplane class; and $t=t_H H$. 
Denote by  $I_{\T}$ the specialization of $I^\zeta_{\T}$
with  $$ \zeta=0. $$
Denote by  $\underline{I}_{\T}$ the specialization of $I_{\T}$ with
\begin{align} \label{root}  \lambda _i = \exp {(2\pi i\sqrt{-1}/n)}, \ \ i=1,..., n. \end{align}

Let's define formal functions $B_k(q,z) \in \QQ[H]/(H^n -1) \otimes_\QQ \QQ[[q,\frac{1}{z}]]$ and $C_k (q)\in \QQ[[q]]$ for $k=0, 1, ..., n-1$ inductively as follows:
First, set $B_0 :=\underline{I}_{\T} (0, q)$ and choose $C_0 (q)$ by a coefficient of $1=H^0$ in $B_0(q,z=\infty)$. 
Now, suppose $B_{k-1}(q,z)$ and $C_{k-1}(q)$ are defined. Then, define $B_k(q,z)$ by 
$$\ B_k   := (H + zq\frac{d}{dq}) \frac{B_{k-1}}{C_{k-1}(q)},$$
and $C_k (q)$ by a coefficient of $H^k$ in $B_k(q,\infty)$. 
We can easily check that $C_k (q)$, $k=0,1, ...,n-1$, which are so called initial constants are of the form $1+O(q)$ and also that 
$$C_k(q)H^k=B_k(q, \infty), ~~~k=0,1,...,n-1.$$ 
Note that $\QQ[H]/(H^n -1)$ is isomorphic to $H^*_\T(\PP^{n-1})$ modulo \eqref{root}.

Now we are ready to state the main theorem of this paper.
\begin{Thm}\label{Main}
\begin{align*}  \lan\; \ran ^{0+}_{1,0}  & =  -\frac{3(n-1-r-m)^2+n-r+m-3}{48} \log (1 - q\prod _{a=1}^r  l_a ^{l_a}\prod _{b=1}^m  (-l'_b) ^{l'_b}) \\
 & - \frac{1}{2}\sum _{k=m}^{n-r-2} \binom{n-r-k}{2}\log C_{k} (q) .  \end{align*}
\end{Thm}

Define  $I^\zeta_0$ and $ I^\zeta_1$ by the $1/z$-expansion of
\[ I^\zeta_\T|_{t = 0} = I^\zeta_0 + I^\zeta_1/z + O(1/ z^2). \]  
Denote $I_0, I_1$ by the specialization of $I^\zeta_0, I^\zeta_1$ with $\zeta=0$. It's easy to check that $I^\zeta_0=I_0=C_0.$
In \cite{CKw}, they proved wall-crossing formula.

\begin{Thm}\cite{CKw} \begin{align*}\label{Wall} 
\lan\; \ran ^{0+}_{1,0}    & = -\frac{1}{24}\chi _{\mathrm{top}} (X) \log I_0 - \frac{1}{24} \int _X \frac{I^\zeta_1}{I_0} c_{\dim X -1} (T_X) \\ & +
   \lan\; \ran ^{\infty}_{1,0}  |_{q^d \mapsto q^d \exp (\int _{d[\mathrm{line}]} \frac{I_1}{I_0} ) } . \nonumber
\end{align*}
\end{Thm}
 
Here, we have to regard $c_{\dim X -1} (T_X)$ as an equivariant Chern class in order to define integration on $X$ by localization.
If we combine those two theorems, we get the following theorem.

\begin{Thm}
\begin{align*}  \lan\; \ran ^{\infty}_{1,0}  |_{q^d\mapsto  q^d \exp (\int _{d[\mathrm{line}]} \frac{I_1}{I_0} )} & = 
\frac{1}{24} \chi _{\mathrm{top}}(X) \log I_0
+ \frac{1}{24}  \int _X \frac{I^\zeta_1}{I_0} c_{\dim X -1} (T_X) 
 \\ & - \frac{3(n-1-r-m)^2+n-r+m-3}{48} \log (1 - q\prod _{a=1}^r  l_a ^{l_a}\prod _{b=1}^r  (-l'_b) ^{l'_b})  
 \\ & - \frac{1}{2}\sum _{k=m}^{n-r-2} \binom{n-r-k}{2}\log C_{k} (q) . \end{align*}
\end{Thm}

When $m=0$, Theorem 1.3 gives another proof of the result in \cite{KL}. Also when $r=0$, it gives another proof of the result in \cite{H}.  If both are non-zero, then this gives a new result.

\subsection{Aknowledgments}

The authors would like to thank Bumsig Kim for initially suggesting the problem, giving many invaluable suggestions, support and advice. 
We also thank Korea Institute for Advanced Study for financial support, excellent working conditions and inspiring research environment.

The research of J. O. was partially supported by the NRF grant 2007-0093859 and by Basic Science Research Program through the National Research Foundation of Korea(NRF) funded by the Ministry of Education(No. NRF-2013R1A1A2007780).

\section{Elliptic quasimap potential function of $X$}

In this section, we will simplify the ellipctic quasimap potential function of $X$.
We closely follow the notations in \cite{KL} and state the results in \cite{KL} without proof in this and next sections.

\subsection{Quantum Lefschetz theorem and Divisor axiom}

We will write elliptic quasi-map potential function as a generating function of quasi-map invariants of $\PP^{n-1}$ which is much easier to deal with.
Consider $Q^{0+}_{g, k}(\PP^{n-1}, d )$. Here and before, since the degree of a quasi-map to $\PP^{n-1}$ can be regarded as a non-negative integer, we used the notation $d$ instead of $\beta$.
Denote by $f$ the universal map
from the universal curve $\mathcal{C}$ of $Q^{0+}_{g, k}(\PP^{n-1}, d )$ to the stack quotient $[\CC^n/\CC^*]$:
\[  \xymatrix{   \mathcal{C} \ar[r]^{f} \ar[d]_{\pi} & [\CC^n/\CC^*] \\
   Q^{0+}_{g, k}(\PP^{n-1}, d )  &    }  .  \]
Since the domain curves of objects in $Q^{0+}_{g, k}(\PP^{n-1}, d )$ don't have rational tails for any $g$ and $k$, every irreducible component with genus $1$ in the domain curves of objects in $Q^{0+}_{1,0}(\PP^{n-1}, d )$ has to have a positive degree if it exists.
So, we can apply Quantum Lefschetz hyperplane section theorem in \cite{Kim} to get following formula: 
$$(j_{g,k,d})_*[Q^{0+}_{g, k} (X, d )]^{\vir} = \re (E_{g,k,d}\oplus E'_{g,k,d}) \cap [Q^{0+}_{g, k} (\PP^{n-1}, d )]^{\vir}$$
for $g=0$, $k=2,3,\cdots$ and $g=1$, $k=0$, $d > 0$. 
Here, 
$$j_{g,k,d} : Q^{0+}_{g, k} (X, d )  \cong  Q^{0+}_{g, k} (X', d )    \hookrightarrow Q^{0+}_{g, k} (\PP^{n-1}, d )$$
if $d>0$ and
$$j_{0,k,0} : Q^{0+}_{0, k} (X, 0 ) \cong X \times \overline{M}_{0,k} \rightarrow \PP^{n-1} \times \overline{M}_{0,k} \cong Q^{0+}_{0, k} (\PP^{n-1}, 0 )$$
where $\overline{M}_{0,k}$ is the moduli stack of stable curves with genus $0$ and $k$-markings.
Also,
$$E_{g,k,d}=R^0\pi _* f^* [(E\times \CC^n)/\CC^*], ~~~~E'_{g,k,d}=R^1\pi _* f^* [(E'\times \CC^{n})/\CC^*]$$ 
where $E=\oplus _a E_a$(resp. $E'=\oplus _b E'_b$), $E_a$(resp. $E'_b$) is one dimensional $\CC^*$-representation space with weight $l_a\theta$ (resp. $-l'_b\theta$) where $\theta:\CC^* \rightarrow \CC^*$ is the identity map. 
And $\re$ stands for the Euler class.
So, we can rewrite the potential function as 
 \begin{align} \label{potF} \lan\; \ran ^{0+ }_{1,0} = \sum_{d \neq 0} q^d \deg \left( \re (E_{1,0,d}\oplus E'_{1,0,d}) \cap [Q^{0+}_{1,0} (\PP^{n-1}, d )]^{\vir} \right). \end{align}



On the other hand, denote by 
$$Q^{0+, 0+}_{g, k|m}(\PP^{n-1}, d)  $$
the moduli space of genus $g$, degree class $d$ stable quasimaps to $\PP^{n-1}$ with
ordinary $k$ pointed markings and infinitesimally weighted $m$ pointed markings.
We can also define their natural virtual fundamental classes
$$[Q^{0+, 0+}_{g, k|m}(\PP^{n-1}, d)]^{vir}$$ 
(see \S 2  of \cite{BigI}). 
Using this, we define invariants
\begin{align*}  & \lan \gamma _1\psi  ^{a_1} , ..., \gamma _k\psi  ^{a_k} ;  \delta _1, ..., \delta _m \ran _{g, k|m, d}^{0+, 0+}  \\ & =
\int _{\re( E_{g,k|m,d}\oplus E'_{g,k|m,d})\cap [Q^{0+, 0+}_{g, k|m} (\PP^{n-1}, d) ]^{\vir}} \prod _i ev_i^*(\gamma _i)\psi _i ^{a_i} \prod _j \hat{ev}_j ^* (\delta _j) \ \ 
\end{align*}
for $\gamma_i \in H^* (\PP^{n-1})\ot \QQ (\lambda )$, $\delta _j \in H^* ([\CC^n/\CC^* ], \QQ )$ and for $g=0$, $k=2,3,\cdots$ and $g=1$, $k=0$, $d > 0$
where $\psi_i$ is the psi-class associated to the $i$-th marking and $ev_i$(resp. $\hat{ev}_j$) is the evaluation map to $\PP^{n-1}$(resp. $[\CC^n/\CC^*]$) at the $i$(resp. $j$)-th marking(resp. infinitesimally weighted  marking);
$E_{g,k|m,d}=R^0\pi _* f^* [(E\times \CC^n)/\CC^*], ~~~~E'_{g,k|m,d}=R^1\pi _* f^* [(E'\times \CC^{n})/\CC^*]$, here, $f$ and $\pi$ are defined as 
\[  \xymatrix{   \mathcal{C} \ar[r]^{f} \ar[d]_{\pi} & [\CC^n/\CC^*] \\
   Q^{0+,0+}_{g,k|m}(\PP^{n-1}, d )  &    }   \]
where $\cC$ is the universal curve.
Note that in the definition of invariants, there is a constraint about $g$ and $k$ because Quantum Lefschetz theorem holds only in this case. So, in this case, we can interpret these invariants as invariants for $X$ which are basically defined without constraint about $g$ and $k$.

In here, we are focusing only on $Q^{0+,0+}_{1,0|1}(\PP^{n-1},d)$ which is isomorphic to the universal curve of $Q^{0+}_{1,0}(\PP^{n-1},d)$. 
Define a generating function
$$\lan  \tilde{H} \ran ^{0+}_{1, 0|1} 
:= \sum _{d=1}^{\infty} q^d \lan ; \tilde{H} \ran ^{0+, 0+}_{1, 0|1, d}$$
where $\tilde{H} \in H^2([\CC^n /\CC^*], \QQ)$ is the hyperplane class.
Then, by divisor axiom, we have
\begin{align} \label{PotF2} q\frac{d}{dq} \lan \; \ran^{0+}_{1,0} = \lan \tilde{H} \ran^{0+}_{1,0|1} \end{align}

\subsection{Localization}

Now, we will calculate it by using $\T$-equivariant quasi-map theory.  
Recall that $\T = (\CC^*)^n$ is $n$-dimensional torus acting on $\PP^{n-1}$ in a standard way.
Let $\{ p_i \}_i$ be the set of $\T$-fixed points of $\PP^{n-1}$.
The $\T$-fixed loci of $Q^{0+,0+}_{1, 0| 1}(\PP^{n-1}, d)$ can be divided into two types according to whether the reduced image is a point in $\PP^{n-1}$ or not.
A quasimap will be called a vertex type over $p_i$ if its regularization map is constant over $p_i$. 
For the definition of regularization map, see \cite{CKM}. Otherwise, the quasimap in $Q^{0+,0+}_{1, 0| 1}(\PP^{n-1}, d)^\T$ will be called a loop type. The loop type quasimap is called a loop type over $p_i$ if the infinitesimally weighted marking of the quasimap maps to $p_i$. 

Define $Q_{vert, i,d} ^\T$ to be the substack of $Q^{0+,0+}_{1, 0|1}(\PP^{n-1}, d) ^\T$ consisting of vertex type over $p_i$.
Define $Q_{loop, i,d} ^\T$ to be the substack of $Q^{0+,0+}_{1, 0|1}(\PP^{n-1}, d) ^\T$ consisting of loop type over $p_i$.

By the virtual localization theorem, 
$ \lan  \tilde{H} \ran ^{0+}_{1, 0|1}  $ can be divided into the sum of the localization contribution
$\mathbf{Vert}_i$ from all the vertex type over $p_i \in {\PP^{n-1}}^\T$ and the localization 
contribution $\mathbf{Loop}_i$
from all the loop type over $p_i\in {\PP^{n-1}}^\T$.
 That is,
\begin{align} \label{PotF3}   \lan  \tilde{H} \ran ^{0+}_{1, 0|1}    & := \sum _i \mathbf{Vert}_i  + \sum _i \mathbf{Loop}_i , \end{align} 
where
\begin{align*}
\mathbf{Vert}_i & := \sum _{d\ne 0} q^d  \int _{[Q_{vert, i,d} ^{\T}]^{\vir}} 
\frac{\ret( E_{1,0|1,d}\oplus  E'_{1,0|1,d})|_{Q_{vert,i,d}^{\T}} \hat{ev}^*_1 (\tilde{H}) }{\ret(N^{vir}_{Q_{vert, i,d}^{\T}/Q^{0+,0+}_{1,0|1}(\PP^{n-1}, d)} ) } ,
\\ \mathbf{Loop}_i  & := \sum_{d\ne 0}  q^d  \int _{[Q_{loop, i,d} ^{\T}]^{\vir}} 
\frac{\ret( E_{1,0|1,d}\oplus  E'_{1,0|1,d})|_{Q_{loop,i,d}^{\T}}  \hat{ev}^*_1 (\tilde{H}) }{\ret(N^{vir}_{Q_{loop, i,d}^{\T}/Q^{0+,0+}_{1,0|1}(\PP^{n-1}, d)} ) }  \end{align*}
where $\ret$ stands for the $\T$-equivariant Euler class and $N^{vir}_{Q_{vert, i,d}^{\T}/Q^{0+,0+}_{1,0|1}(\PP^{n-1}, d)}$(resp. $N^{vir}_{Q_{loop, i,d}^{\T}/Q^{0+,0+}_{1,0|1}(\PP^{n-1}, d)}$) is the virtual normal bundle of $Q_{vert, i,d}^{\T}$(resp. $Q_{loop, i,d}^{\T}$) into $Q^{0+,0+}_{1,0|1}(\PP^{n-1}, d)$.
Here, we regard the hyperplane class $\tilde{H}$ as a $\T$-equivariant class in $H^2_{\T}([\CC^n /\CC^*], \QQ).$

\subsubsection{Vertex term}

Let 
$$Q^{0+}_{1,0}(\PP^{n-1},d)^{\T,p_i}$$
be the $\T$-fixed part of $Q^{0+}_{1,0}(\PP^{n-1},d)$ whose elements have domain components only over $p_i$ under the regularization map.
Then, $Q^{\T}_{vert,i,d}$ is isomorphic to the universal curve of $Q^{0+}_{1,0}(\PP^{n-1},d)^{\T,p_i}$.
So, by the divisor axiom,
$$\mathbf{Vert}_i =q\frac{d}{dq} \sum_{d\ne 0} q^d \int_{[Q^{0+}_{1,0}(\PP^{n-1},d)^{\T,p_i}]^{vir}}  
\frac{\ret( E_{1,0,d}\oplus  E'_{1,0,d})|_{Q^{0+}_{1,0}(\PP^{n-1},d)^{\T,p_i}}  }{\ret(N^{vir}_{Q^{0+}_{1,0}(\PP^{n-1},d)^{\T,p_i}/Q^{0+}_{1,0}(\PP^{n-1},d)} ) }$$
By the way,
\begin{align*} \ret( E_{1,0,d})|_{Q^{0+}_{1,0}(\PP^{n-1},d)^{\T,p_i}} & = \ret(\prod _a \pi_* \mathcal{O}_{ \cC}(l_a\hat{\bx}) \ot E_a ) \\
 & = \ret(\prod _a \pi_* \mathcal{O}_{ l_a\hat{\bx}}(l_a\hat{\bx}) \ot E_a) \ret(\prod_a R\pi_* \cO_{\cC} \otimes E_a)\\
 &=  \ret(\prod _a \pi_* \mathcal{O}_{ l_a\hat{\bx}}(l_a\hat{\bx}) \ot E_a) \frac{\ret(\prod_a \pi_* \cO_{\cC} \otimes E_a)}{\ret(\prod_a R^1\pi_* \cO_{\cC} \otimes E_a)}\end{align*}
where $\hat{\bx}$ is base loci on a universal curve $\cC$ and $\pi$ is a projection from $\cC$ to $Q^{0+}_{1,0}(\PP^{n-1},d)^{\T,p_i}$.
First equatlity comes from the idea in \cite{CKg}. And second equality comes from the long exact sequence
$$0 \ra \cO_\cC \ra \mathcal{O}_{ \cC}(l_a\hat{\bx}) \ra \mathcal{O}_{ l_a\hat{\bx}}(l_a\hat{\bx})  \ra 0.$$
Similarly, we can show that 
\begin{align*} \ret( E'_{1,0,d})|_{Q^{0+}_{1,0}(\PP^{n-1},d)^{\T,p_i}}  = \ret(\prod _b \pi_* \mathcal{O}_{ l'_b\hat{\bx}} \ot E'_b)\frac{\ret(\prod_b R^1\pi_* \cO_{\cC} \otimes E'_b)}{ \ret(\prod_b \pi_* \cO_{\cC} \otimes E'_b)}. \end{align*}
And also, we can see that
\begin{align*} \ret(N^{vir}_{Q^{0+}_{1,0}(\PP^{n-1},d)^{\T,p_i}/Q^{0+}_{1,0}(\PP^{n-1},d)} )& = \ret(R\pi_*  \mathcal{O}_{\cC}(\hat{\bx} ) \ot T_{p_i} \PP^{n-1}) \\
 & = \ret(\pi_* \mathcal{O}_{\hat{\bx}}(\hat{\bx} ) \ot T_{p_i} \PP^{n-1})\frac{ \ret( \pi_* \cO_{\cC} \otimes T_{p_i} \PP^{n-1})}{\ret( R^1\pi_* \cO_{\cC} \otimes T_{p_i} \PP^{n-1})}. \end{align*}
On the other hand, 
$$Q^{0+}_{1,0}(\PP^{n-1},d)^{\T,p_i} \cong \overline{M}_{1, 0|d} /S_d$$
where $S_d$ is the symmetric group of degree $d$ acting on $\overline{M}_{1, 0|d}$ by a permutation of infinitesimally weighted markings.
Furthermore, $\overline{M}_{1, 0|d}$ is smooth and
$$[Q^{0+}_{1,0}(\PP^{n-1},d)^{\T,p_i}]^{vir} = \frac{1}{d!}[\overline{M}_{1, 0|d}]$$
under the isomorphism. Here, $[\overline{M}_{1, 0|d}]$ is the fundamental class of $\overline{M}_{1, 0|d}$.
Therefore,
$$
         \mathbf{Vert}_i  = q\frac{d}{dq}\sum _{d \ne 0} \frac{q^d}{d!}  \int _{\overline{M}_{1, 0|d}} (1+ c_i(\lambda)    \re (\mathbb{E}) ) F_{i,d} $$
where
\begin{align} \label{Fterm}
F_{i,d}  :=
\frac{\ret(\prod _a \pi_* \mathcal{O}_{ l_a\hat{\bx}}(l_a\hat{\bx}) \ot E_a) \ret(\prod _b \pi_* \mathcal{O}_{ l'_b\hat{\bx}} \ot E'_b) }{\ret(\pi_* \mathcal{O}_{\hat{\bx}}(\hat{\bx} ) \ot T_{p_i} \PP^{n-1}) } \end{align}
with $\hat{\bx} := \sum_{j=1}^d \hat{x}_j$, sum of loci of infinitesimally weighted markings in the universal curve; $\mathbb{E} := (R^1\pi_* \cO_{\cC})^\vee$ is the Hodge bundle on $\overline{M}_{1, 0|d}$ and $c_i(\lambda)$ is the element in  $\QQ (\lambda)$ uniquely determined by 
\begin{align*} 1+ c_i(\lambda) \re (\mathbb{E})  
&  = \frac{\re ^\T(\mathbb{E}^\vee \ot T_{p_i}\PP^{n-1} ) \re ^\T ( \cO_{\overline{M}_{1, 0|d}} \otimes E |_{p_i}) \re ^\T (\mathbb{E}^\vee \ot E'|_{p_i} ) }
{\re ^\T(\cO_{\overline{M}_{1, 0|d}} \otimes T_{p_i} \PP^{n-1}) \re ^\T (\mathbb{E}^\vee \ot E|_{p_i} ) \re ^\T ( \cO_{\overline{M}_{1, 0|d}} \otimes E' |_{p_i})}. 
\end{align*}
So, by a simple computation, we can see that
$$c_i(\lambda)  =  \sum_{j\ne i} \frac{1}{\lambda _j - \lambda _i} + \sum _a \frac{1}{l_a\lambda _i}+\sum _b \frac{1}{l'_b\lambda _i}.$$
Note that $c_i(\lambda)$ is independent of $d$ and that
$\re(\mathbb{E})^2 =0$ because $\mathbb{E}$ comes from the Hodge bundle on $\overline{M}_{1,1}$.

By the same argument in \cite{KL}, we can relate the genus one invariants with the genus zero invariants. 

\begin{Prop}
$$ 24 \sum_{d \neq 0} \frac{q^d}{d!} \int _{\overline{M}_{1, 0|d}} \re (\mathbb{E})  F_{i, d} 
       = \sum_{d \neq 0} \frac{q^d}{d!} \int _{\overline{M}_{0, 2| d}} F_{i, d},  $$

$$e^{24\sum_{d\ne 0}  \frac{q^d}{d!} \int _{\overline{M}_{1, 0|d}} F_{i, d}}
      = \sum_{d \neq 0} \frac{q^d}{d!} \int _{\overline{M}_{0, 3| d}} F_{i, d}  $$
where $F_{i,d}$ classes on $\overline{M}_{0,2|d}$ and $\overline{M}_{0,3|d}$ are defined in the same way as \eqref{Fterm} 
\end{Prop}

In conclusion, we have
\begin{align} \label{VerF}
         \mathbf{Vert}_i  = \frac{q}{24} \frac{d}{dq} \left( c_i(\lambda)\sum _{d \ne 0} \frac{q^d }{d!}  \int _{\overline{M}_{0,2|d}}  F_{i,d} ~~ + \text{log}\left( \sum _{d \ne 0} \frac{q^d}{d!}  \int _{\overline{M}_{0,3|d}}  F_{i,d} \right) \right). \end{align}

\section{Localized invariants}

\subsection{Localized generating functions in genus zero theory}

In order to do equivariant quasi-map theory of $\PP^{n-1}$ instead of that of $X$, we need to use $E \times E'$-twisted Poincar\'e metric on $H^*_\T(\PP^{n-1})\ot \QQ (\lambda)$. That is for $a,b \in H^*_\T(\PP^{n-1})\ot \QQ (\lambda)$,
$$\lan a,b \ran^{E \times E'} = \int _{\PP^{n-1}} \frac{\ret (E) \cup a \cup b} {\ret (E')} $$
where $\ret (E)$(resp. $\ret (E')$) is the $\T$-equivariant Euler class of $E$(resp. $E'$).
Here, we used the notation $E$(resp. $E'$) instead of 
$[(E \times (\CC^n \backslash \{0\})) / \CC^* ] \cong \oplus_a \cO_{\PP^{n-1}}(l_a)$(resp. $[(E' \times (\CC^n \backslash \{0\})) / \CC^* ] \cong \oplus_b \cO_{\PP^{n-1}}(-l'_b)$) by avoiding abuse of notations. 
Let $\phi_i$ be the basis of $H^*_\T(\PP^{n-1})\ot \QQ (\lambda)$ such that
$$\phi_i|_{p_j}= \left\{ \begin{array}{rl} 1 & \text{if } i=j  \\
                                    0 & \text{if } i\ne j   \, ,\end{array}\right.$$
and let $\phi^i$ be its dual basis with respect to $E \times E'$-twisted Poincar\'e metric.

Also as in \cite{KL}, we need to use the twisted virtual fundamental class
$$\ret (E_{0,k,d} \oplus E'_{0,k,d})\cap [Q^{0+}_{0, k}(\PP^{n-1}, d)]^{\vir}$$
in genus zero quasi-map theory. By using this, we will define local correlators.
Let $$Q^{0+}_{0, k} (\PP^{n-1}, d) ^{\T, p_i}$$  be the $\T$-fixed part of $Q^{0+}_{0, k} (\PP^{n-1}, d) $ 
whose elements have domain components only over $p_i$. 
By using the twisted virtual fundamental class
\[ \frac{\ret (E_{0,k,d}\oplus E'_{0,k,d}) \cap [Q^{0+}_{0, k} (\PP^{n-1}, d) ^{\T, p_i}]^{\vir}}{\ret (N^{\vir}_{Q^{0+}_{0, k} (\PP^{n-1}, d) ^{\T, p_i}/
Q^{0+}_{0, k}(\PP^{n-1}, d)} )}, \]
define it as follows:
\begin{align*}  & \lan \gamma _1\psi  ^{a_1} , ..., \gamma _k\psi  ^{a_k} \ran _{0, k, d}^{0+, p_i}  :=
\int _{\frac{\ret(E_{0,k,d}\oplus E'_{0,k,d})\cap [Q^{0+}_{0, k} (\PP^{n-1}, d) ^{\T, p_i}]^{\vir}}{\ret (N^{\vir}_{Q^{0+}_{0, k} (\PP^{n-1}, d) ^{\T, p_i}/
Q^{0+}_{0, k}(\PP^{n-1}, d)} )}} \prod _i ev_i^*(\gamma _i)\psi _i ^{a_i} \ \ ; \\
&  \lla \gamma _1\psi  ^{a_1} , ..., \gamma _k\psi  ^{a_k} \rra _{0, k}^{0+, p_i} \\
& := \sum _{m, d} \frac{q^{d}}{m!}
 \lan    \gamma _1\psi  ^{a_1} , ..., \gamma _k\psi  ^{a_k} , t, ..., t  \ran_{0, k+m, d}^{0+, p_i} , \text{ for } t \in H_{\T}^* (\PP^{n-1})\ot \QQ (\lambda ) \ ,
      \end{align*} 
      where $\psi _i$ is the psi-class associated to the $i$-th marking; $ev_i$ is the 
$i$-th evaluation map and $q$ is a formal Novikov variable.
Here, we used the notation $E_{0,k,d}$(resp. $E'_{0,k,d}$) instead of 
$E_{0,k,d}|_{Q^{0+}_{0, k} (\PP^{n-1}, d) ^{\T, p_i}}$(resp. $E'_{0,k,d}|_{Q^{0+}_{0, k} (\PP^{n-1}, d) ^{\T, p_i}} $) by avoiding abuse of notations.

Let $z$ be a formal variable.
We define the following $\T$-local generating functions:
\begin{align*} D_i & : = e_i  \lla 1, 1, 1 \rra ^{0+, p_i}_{0, 3 } = 1 + O(q) \ \ ; \\
  u_i & :=  e_i \lla  1, 1 \rra_{0, 2}^{0+, p_i}    =  t|_{p_i} + O(q)  \ \ ; \\
 S_t^{0+, p_i} (\gamma ) & := e_i
 \lla \frac{1}{z-\psi } , \gamma \rra_{0, 2}^{0+, p_i} = e^t\gamma |_{p_i} + O(q) \\ 
        & \text{ for } \gamma \in H^*_{\T}(\PP^{n-1})\ot \QQ (\lambda ) [[q]]  \ \ ; \\
 J^{0+, p_i} &:= e_i \lla \frac{1 }{z(z-\psi ) }\rra _{0,1}^{0+, p_i} =  e^t |_{p_i} + O(q) \ ,
\end{align*} 
where the unstable terms of $S_t^{0+, p_i} (\gamma )$ and $J^{0+, p_i}$ are defined by using the quasimap graph spaces $QG^{0+}_{0, 0, d} (\PP^{n-1})$ or $QG^{0+}_{0, 1, 0} (\PP^{n-1})$ as in \cite{CK, CKg0}.
Also, the unstable term of $u_i$ (this is the only case of $m=d=0$) is defined to be $0$.
So, in particular,
\begin{align} \label{locI} 
J^{0+, p_i}|_{t=0}  = I_{\T} |_{t=0,\, p_i}.
\end{align}
Here the front terms $e_i$ are defined by the formulas $\phi^i=e_i \phi_i$.
The parameter $z$ naturally appears as the $\CC ^*$-equivariant parameter in the graph 
construction (see \S 4 of \cite{CKg0}). It is originated from the $\CC ^*$-action on $\PP ^1$.

By the way, 
it is easy to check that
$$ \sum_{d \neq 0} \frac{q^d}{d!} \int _{\overline{M}_{0, 2| d}} F_{i, d} = u_i |_{t=0} 
\text{,    }
 \sum_{d \neq 0} \frac{q^d}{d!} \int _{\overline{M}_{0, 3| d}} F_{i, d} = D_i |_{t=0}.$$
So, applying it to \eqref{VerF}, we have
\begin{align}  \label{VerF2}  \mathbf{Vert}_i  = \frac{q}{24} \frac{d}{dq} \left( \left( \sum_{j\ne i} \frac{1}{\lambda _j - \lambda _i} + \sum _a \frac{1}{l_a\lambda _i}+\sum _b \frac{1}{l'_b\lambda _i} \right)u_i |_{t=0} + \text{log} D_i |_{t=0} \right). \end{align}

In order to describe this by using $I$-function for $X$, we need more generating functions. Denote by $QG^{0+}_{0, k, d } (\PP^{n-1}) $ the quasimap graph spaces (see \cite{CKg0})
and by $$QG^{0+}_{0, k, d } (\PP^{n-1})^{\T, p_i} $$ the $\T$-fixed part of $QG^{0+}_{0, k, d } (\PP^{n-1}) $ 
whose elements have domain components only over $p_i$. 
As in \cite{KL}, we define invariants and generating functions on the graph spaces:
for $\gamma _i \in H^*_{\T}(\PP^{n-1})\ot H^*_{\CC^*} (\PP ^1) \ot \QQ (\lambda )$
\begin{align*}
& \lan \gamma _1\psi  ^{a_1} , ..., \gamma _k\psi  ^{a_k} \ran _{k, d}^{QG^{0+}, p_i} =
\int _{\frac{\ret(R^0\pi_*f^* E\oplus R^1\pi_*f^* E')\cap [QG^{0+}_{0, k, d} (\PP^{n-1}) ^{\T, p_i}]^{\vir}}{\ret (N^{\vir}_{QG^{0+}_{0, k, d} (\PP^{n-1}) ^{\T, p_i}/
QG^{0+}_{0, k, d }(\PP^{n-1})} )}} \prod _i ev_i^*(\gamma _i)\psi _i ^{a_i} \ \ ; \\
&  \lla \gamma _1\psi  ^{a_1} , ..., \gamma _k\psi  ^{a_k} \rra _{k}^{QG^{0+}, p_i} \\
& = \sum _{m, d} \frac{q^{d}}{m!}
 \lan    \gamma _1\psi  ^{a_1} , ..., \gamma _k\psi  ^{a_k} , t, ..., t  \ran_{k+m, d}^{QG^{0+}, p_i} , \text{ for } t \in H_{\T}^* (\PP^{n-1})\ot \QQ (\lambda ) .
\end{align*}
Here we denote by $ev_i$ the $i$-th evaluation map to $\PP^{n-1}\times \PP ^1$ from the quasimap graph spaces
and regard $t$ as the element $t\ot 1$  in $H^*_{\T}(\PP^{n-1})\ot H^*_{\CC^*} (\PP ^1) \ot \QQ (\lambda )$.
We used the notation $E$(resp. $E'$) instead of 
$[(E \times (\CC^n )) / \CC^* ] $(resp. $[(E' \times (\CC^n )) / \CC^* ] $) by avoiding abuse of notations.
$f$ and $\pi$ are defined as follows:
\[ \xymatrix{   \mathcal{C} \ar[r]^{f} \ar[d]_{\pi} & [\CC^n/\CC^*] \\
   QG^{0+}_{0, k,d}(\PP^{n-1} )^{\T,p_i}  &    }  .\]
Here, $\cC$ is the universal curve.

Let ${\bf p}_\infty$ be the equivariant cohomology class in
$H^*_{\CC ^*} (\PP ^1)$ defined by 
\[ {\bf p}_{\infty} |_{0} = 0, {\bf p}_{\infty} |_{\infty} = - z .\]
Exactly same as in \cite{KL}, we can have the following factorization.

\begin{Prop} \label{fact}
\begin{equation*}\label{Der_Dr} J^{0+, p_i} = S_t ^{0+, p_i} (P^{0+, p_i}) , \end{equation*}
where 
\[ P^{0+, p_i} :=   e_i  \lla 1 \ot {\bf p}_\infty \rra _{1 }^{QG^{0+}, p_i}. \]
\end{Prop}

By the uniqueness lemma in \S 7.7 of \cite{CKg0}, 
\begin{equation*}\label{localS}  S_t^{0+, p_i} (\gamma ) = e^{u_i/z} \gamma |_{p_i} . \end{equation*}
Hence Proposition \ref{fact} gives
the expression 
\begin{align} \label{Jftn} J^{0+, p_i} = e^{u_i/z} (r_{i, 0} + O(z)), \end{align} where 
$r_{i, 0}\in \QQ (\lambda ) [[t, q]]$ is the constant term of $P^{0+, p_i}$ in $z$.
By the following result, we can easily see that the expression \eqref{Jftn} is unique.

\begin{Cor} The equality 
\[ \log J^{0 +, p_i} = u_i/z + \log r_{i, 0} + O(z)  \ \ \in \QQ (\lambda ) ((z)) [[t, q]] \]
holds  as  Laurent series of $z$ over the
coefficient ring $\QQ (\lambda)$ in each power expansion of $t$ and $q$, after 
regarding $t$ as a formal element.
\end{Cor}

Also as in \cite{KL}, we can have the following result again.

\begin{Cor}
\begin{align*}\label{Dr} D_i |_{t=0}= \frac{1}{r_{i, 0}|_{t=0}}  . \end{align*}
\end{Cor}

In conclusion, applying these to \eqref{VerF2}, we have
\begin{equation} \label{VerF3}   \mathbf{Vert}_i  = \frac{q}{24} \frac{d}{dq} \left( \left( \sum_{j\ne i} \frac{1}{\lambda _j - \lambda _i} + \sum _a \frac{1}{l_a\lambda _i}+\sum _b \frac{1}{l'_b\lambda _i} \right)u_i |_{t=0} - \text{log} (r_{i,0} |_{t=0}) \right) \end{equation}
where $u_i$ and $r_{i,0}$ are defined in terms of factors in $J^{0+,p_i}$ as in \eqref{Jftn}.
Also, $u_i |_{t=0}$ and $r_{i,0} |_{t=0}$ are related to factors in $I_\T$ by \eqref{locI}.

To describe vertex term more concretely, we need further.
Denote by $$ Q^{0+, 0+}_{g, k|m}(\PP^{n-1}, d )^{\T, p_i}$$
the $\T$-fixed 
part of $Q^{0+, 0+}_{g, k|m}(\PP^{n-1}, d)$ whose domain components are
only over $p_i$.

For $\gamma_i \in H^*_{\T} (\PP^{n-1})\ot \QQ (\lambda )$, $\tilde{t}, \delta _j \in H^*_{\T} ([\CC^n/\CC^* ], \QQ )$ denote

\begin{align*}  & \lan \gamma _1\psi  ^{a_1} , ..., \gamma _k\psi  ^{a_k} ;  \delta _1, ..., \delta _m \ran _{0, k|m, d}^{0+, 0+}  \\ & =
\int _{\ret(E_{0,k|m,d}\oplus E'_{0,k|m,d})\cap [Q^{0+, 0+}_{0, k|m} (\PP^{n-1}, d) ]^{\vir}} \prod _i ev_i^*(\gamma _i)\psi _i ^{a_i} \prod _j \hat{ev}_j ^* (\delta _j) \ \ ; \\
&  \lla \gamma _1\psi  ^{a_1} , ..., \gamma _k\psi  ^{a_k} \rra _{0, k}^{0+, 0+} \\
& = \sum _{m, d} \frac{q^{d}}{m!}
 \lan    \gamma _1\psi  ^{a_1} , ..., \gamma _k\psi  ^{a_k} ; \tilde{t}, ..., \tilde{t}  \ran_{0, k|m, d}^{0+, 0+} \ \  ; \\ 
 & \lan \gamma _1\psi  ^{a_1} , ..., \gamma _k\psi  ^{a_k} ;  \delta _1, ..., \delta _m \ran _{0, k|m, d}^{0+, 0+, p_i}  \\ & =
\int _{\frac{\ret( E_{0,k|m,d}\oplus  E'_{0,k|m,d})\cap [Q^{0+, 0+}_{0, k|m} (\PP^{n-1}, d) ^{\T, p_i}]^{\vir}}{\ret (N^{\vir}_{Q^{0+, 0+}_{0, k|m} (\PP^{n-1}, d) ^{\T, p_i}/
Q^{0+, 0+}_{0, k|m}(\PP^{n-1}, d)} )}} \prod _i ev_i^*(\gamma _i)\psi _i ^{a_i} \prod _j \hat{ev}_j ^* (\delta _j) \ \  ; \\
&  \lla \gamma _1\psi  ^{a_1} , ..., \gamma _k\psi  ^{a_k} \rra _{0, k}^{0+, 0+, p_i} \\
& = \sum _{m, d} \frac{q^{d}}{m!}
 \lan    \gamma _1\psi  ^{a_1} , ..., \gamma _k\psi  ^{a_k} ; \tilde{t}, ..., \tilde{t}  \ran_{0, k|m, d}^{0+, 0+, p_i}. 
\end{align*}

Consider
\begin{align*}
\mathds{S}(\gamma) & : =  \sum _{i}  \phi ^i  \lla \frac{\phi _i}{z-\psi} , \gamma \rra _{0, 2}^{0+, 0+}\ \ ;  \\
\mathds{V}_{ii} (x, y)  & := 
\lla \frac{\phi_i}{x- \psi } ,  \frac{\phi_i}{y - \psi } \rra _{0, 2}^{0+, 0+} 
= \frac{1}{e_i(x+y)} + O(q) \ \ ; \\
 \mathds{U}_i & :=  e_i 
 \lla 1, 1 \rra_{0, 2}^{0+, 0+,  p_i}    =  \tilde{t}|_{p_i} + O(q)  \ \ ; \\
\mathds{S}_i^{0+, p_i} (\gamma) & := e_i 
 \lla \frac{1}{z-\psi} , \gamma \rra _{0, 2}^{0+, 0+,  p_i}   = e^{\tilde{t}/z}\gamma |_{p_i}+ O(q)\ \ ;\\
\mathds{J}^{0+, p_i} &:= e_i 
 \lla \frac{1}{z(z-\psi)}  \rra_{0, 1 }^{0+, 0+, p_i} = e^{\tilde{t}}|_{p_i} + O(q) = J^{0+, p_i}|_{t=0} + O(\tilde{t}) .
 \end{align*}

As before, 
\begin{align} \mathds{S}_i^{0+, p_i} (\gamma)  &= e^{\mathds{U}_i/z} \gamma |_{p_i} \ \ ; \nonumber \\
\label{asymJ}  \mathds{J}^{0+, p_i} & = e^{\mathds{U}_i/z} (\sum_{k=0}^m \mathds{R}_{i,k}z^k + O(z^{m+1})) 
\end{align} 
for some $\mathds{R}_{i,k} \in \QQ (\lambda)[[\tilde{t}, q]]$ (after regarding $\tilde{t}$ as a formal 
element).

Denote by $\mathds{I}$ the infinitesimal $I$-function $\mathds{J}^{0+, 0+}$ defined and calculated explicitly in \cite{BigI}:

\begin{equation*}\mathds{I} (\tilde{t} ) =  \left(\exp (\sum _{i=0}^{n-1} \frac{t_i}{z} (z q \frac{d}{dq} + H)^i ) \right) I_{\T} |_{t=0},
\end{equation*}
By \eqref{asymJ} and the fact that $\mathds{I} |_{ p_i}=\mathds{J}^{0+,p_i}$, we have
$$\mathds{I}|_{p_i} = e^{\mathds{U}_i/z} (\sum_{k=0}^m \mathds{R}_{i,k}z^k + O(z^{m+1})). $$
Hence,
$$\mathds{I}|_{\tilde{t}=0,p_i} = e^{\mathds{U}_i|_{\tilde{t}=0}/z} (\sum_{k=0}^m \mathds{R}_{i,k}|_{\tilde{t}=0}z^k + O(z^{m+1})). $$
On the other hand, since $\mathds{I} |_{\tilde{t}=0} = I_{\T} |_{t=0} \in H^*_\T (\PP^{n-1})[[q,\frac{1}{z}]]$ and $I_{\T} |_{t=0}$ is homogeneous of degree $0$ if we set
$$\text{deg}H=\text{deg}\lambda_i =\text{deg}z =1~~ \text{ and }~~ \text{deg}q=0 ,$$ 
after the specialization
 \begin{align} 
 \label{spec}
 \lambda _i = \lambda_0 \exp {(2\pi i\sqrt{-1}/n)}, \ \ i=1,..., n, \ \ \text{deg}\lambda_0 =1, \end{align}
$\mathds{I} |_{\tilde{t}=0} \in \QQ[H]/(H^n-\lambda_0^n)[[q,\frac{1}{z}]]$ is also homogeneous of degree $0$. 
Since $\mathds{I} |_{\tilde{t}=0, p_i} = \mathds{I} |_{\tilde{t}=0, H=\lambda_i}$ and $\lambda^n_i=\lambda_0^n$, $\mathds{I} |_{\tilde{t}=0, p_i}$ modulo (\ref{spec}) is a series in $\QQ[[q,\frac{\lambda_i}{z}, \frac{z}{\lambda}]]$. 
So, we obtain 
\begin{align}  
\label{gft}
 \mathds{I} |_{\tilde{t}=0, p_i} 
\equiv   e^{\mu (q)\lambda _i /z}  (\sum _{k=0}^{\infty} R_{k}(q) (z /\lambda _i)^k )  \end{align}
for some  $\mu (q) \in \QQ [[q]] $ and  $R_{k}(q) \in \QQ [[q]]$.
Here, $\equiv$ means modulo (\ref{spec}).
Even if we put $\lambda_0 =1$ to match the specialization (\ref{spec}) with (\ref{root}), (\ref{gft}) still holds after (\ref{root}).
Since $\mu(q)\lambda_i = \mathds{U}_i |_{\tilde{t}=0}$ modulo (\ref{spec}), $\mu(q) \in q \QQ[[q]].$
Note that $\mu(q)$ and $R_k(q)$ are independent of $i$.
Hence
\begin{align*}   \mathds{I} |_{\tilde{t}=t_H\tilde{H}, p_i} & = e^{\kl_it_H /z}  (I_{\T} |_{t=0, q\mapsto qe^{t_H},p_i} ) \\
& \equiv  e^{\kl_it_H /z}  e^{\mu (qe^{t_H})\lambda _i /z}  (\sum _{k=0}^{\infty} R_{k}(qe^{t_H}) (z /\lambda _i)^k ) .   \end{align*}
Thus
\begin{align} \label{UandR}  \mathds{U}_i |_{\tilde{t}=t_H  \tilde{H} } & \equiv  \lambda _i (t_H + \mu (qe^{t_H})) \text{ and }  \\
      \mathds{R}_{i, k}|_{\tilde{t}=t_H  \tilde{H} } & \equiv  R_{k}(qe^{t_H})/ (\lambda _i)^k . \nonumber \end{align}
Since 
\begin{align*} r_{i, 0}|_{t=0} & = \mathds{R}_{i, 0} |_{\tilde{t}=0} , \\ 
u_i|_{t=0} &  = \mathds{U}_i|_{\tilde{t}=0} , \text { and } \\
c_i(\lambda) & =  (\sum_{j\ne i} \frac{1}{\lambda _j - \lambda _i} )+ \sum _a \frac{1}{l_a\lambda _i}+\sum _b \frac{1}{l'_b\lambda _i}, \end{align*}
we conclude that
\begin{align}  \label{Vert}
&\sum_i \mathbf{Vert}_i = q\frac{q}{dq}\left(\sum _i \frac{-\log r_{i,0}|_{t=0}}{24} + \sum _i\frac{c_i(\lambda) u_i |_{t=0}}{24}\right) \\
                                    & \equiv q\frac{q}{dq}\left(\frac{ -n \log  R_0(q)}{24} + \frac{1}{24} (\sum_a  \frac{n}{l_a}+\sum_b  \frac{n}{l'_b} -\binom{n}{2}) \mu (q) \right) \nonumber \end{align}
where $\mu$ and $R_0$ are defined in terms of factors in $I_\T$.

\subsection{Loop term}

By the same argument in \cite{Elliptic}, we can represent the loop term with genus zero invariants.
\begin{Prop}
\begin{align} \label{LoopF} \mathbf{Loop}_i & = \frac{1}{2} ( \frac{d}{dt_H} \mathds{U}_i |_{\tilde{t}=t_H \tilde{H}}) |_{t_H=0} 
\lim _{x, y\mapsto 0} \left(e^{-\mathds{U}_i(\frac{1}{x} + \frac{1}{y})} e_i\mathds{V}_{ii}(x, y) - \frac{1}{x+y} \right)  |_{\tilde{t}=0}\\
& \equiv \frac{1}{2}  \lambda_i (1+ q\mu'(q))
\lim _{x, y\mapsto 0} \left(e^{-\lambda_i \mu(q)(\frac{1}{x} + \frac{1}{y})} e_i\mathds{V}_{ii}(x, y)|_{\tilde{t}=0} - \frac{1}{x+y} \right). \nonumber  \end{align} 
\end{Prop}

The second equivalence comes from \eqref{UandR}.

Now consider an equivariant cohomology basis $$\{ 1, H:=c_1^{\T}(\mathcal{O}(1)), ..., H^{n-1}\}$$ of the $\T$-equivariant cohomology ring
\begin{align*} H^*_\T(\PP ^{n-1}) \cong \QQ [\lambda _1, ..., \lambda _n, h]/ (\prod _{i=1}^n  (h - \lambda _i) ) , \ \ \ 
                           H \mapsto  h . \end{align*}

There is the expression of $V$-correlators in terms of $S$-correlators by Theorem 3.2.1 of \cite{CKg}:
\begin{align*}  e_i \mathds{V}_{ii} (x, y)|_{\tilde{t}=0}  = \frac{1}{e_i}
\frac{\sum _j \mathds{S}_{z=x, \tilde{t}=0} (H^j)|_{p_i} \mathds{S}_{z=y, \tilde{t}=0} ({H^j}^\vee )|_{p_i}}{x+ y}  , \end{align*}
where ${H^j}^\vee$ is the daul of $H^j$ with respect to $E \ti E'$-twisted Poincar\'e metric $g_{ij}$ modulo relations \eqref{spec};

\[ g_{ij} = \frac{\prod _{a=1}^r (l_a) }{\prod _{b=1}^m (-l'_b)} \sum _{k=0}^{2} \lambda_0^{n(k-1)} \delta_{r-m+i+j, nk-1} 
\text{ for } 0\le i, j \le n-1.\]

We can calculate the $\mathds{S}_{\tilde{t}=0}(H^k)$ in terms of $I_\T$ by Birkhoff factorization method as in \cite{KL}.
So, we can express $e_i \mathds{V}_{ii} (x,y)|_{\tilde{t}=0}$ and $\mathbf{Loop}_i$ in terms of factors of $I_\T$.

By using Birkhoff factorization method and the calculation in \cite{Popa1}, we obtain
\begin{align*} \lim _{x, y\mapsto 0} \left(e^{-\lambda_i \mu(q)(\frac{1}{x} + \frac{1}{y})} e_i\mathds{V}_{ii}(x, y)|_{\tilde{t}=0} - \frac{1}{x+y} \right)=\frac{\lambda^{n-2}_i}{\ret(T_{p_i} \PP^{n-1})L(q)}q\frac{d}{dq}\text{Loop}(q) \end{align*}
where
\begin{align*} 
& L(q)=(1-q\prod_a l_a^{l_a}\prod_b (-l'_b)^{l'_b})^{-\frac{1}{n}}, \\
& \text{Loop}(q) = 
  \frac{n}{24} ( n-1 - 2 \sum_{a=1}^r \frac{1}{l_a} - 2 \sum_{b=1}^m \frac{1}{l'_b}) \mu (q)    
   \\ & - \frac{3(n-1-r-m)^2 + (n-2)}{24}\log (1-q\prod_a l_a^{l_a}\prod_b (-l'_b)^{l'_b})  
   \\ & - \sum _{k=m}^{n-r-2}\binom{n-r-k}{2} \log C_k (q) . \end{align*}
Thus,
\begin{align} \label{LoopF} \mathbf{Loop}_i   \equiv \frac{1}{2}  \lambda_i (1+ q\mu'(q))
\frac{\lambda^{n-2}_i}{\ret(T_{p_i} \PP^{n-1})L(q)}q\frac{d}{dq}\text{Loop}(q).   \end{align}

Using the fact that equivariant $I_{\T}$-function satisfies the Picard-Fuchs equation
\[ \mathrm{PF} I_{\T} |_{t=t_H \cdot H} =0, \ \ \mathrm{PF}:= ( z\frac{d}{dt} )^n - 1 - 
q\prod _{a} \prod _{m=1}^{l_a} (l_az\frac{d}{dt} + mz) \prod _{b} \prod _{m=0}^{l'_b-1} (-l'_bz\frac{d}{dt} - mz) , \]
and asymptotic form of \[ I_{\T} |_{t=t_H H, p_i},  \]
we can calculate $\mu$ and $R_0$, 
\begin{align*}
\mu(q) & = \int_0^q \frac{ L(x)- 1}{x} dx,  ~~R_0(q) = L(q)^{\frac{r-m+1}{2}}.
\end{align*}
For calculations, see \cite{Popa1}. Therefore,
\begin{align}
\label{Loop1}\sum_i \mathbf{Loop}_i & \equiv \frac{1}{2} q\frac{d}{dq}\text{Loop}(q) \sum_i \frac{\lambda_i^{n-1}}{\ret(T_{p_i} \PP^{n-1})}\\
& = \frac{1}{2} q\frac{d}{dq}\text{Loop}(q) \int_{\PP^{n-1}} H^{n-1} = \frac{1}{2} q\frac{d}{dq}\text{Loop}(q). \nonumber
\end{align}

\subsection{Proof of main theorem}

By combining \eqref{Vert}, \eqref{Loop1} with \eqref{PotF2}, we have
\begin{align*}  \frac{d}{dq} \{ \lan\; \ran ^{0+}_{1,0}  & +\frac{3(n-1-r-m)^2+n-r+m-3}{48} \log (1 - q\prod _{a=1}^r  l_a ^{l_a}\prod _{b=1}^m  (-l'_b) ^{l'_b}) \\
 & + \frac{1}{2}\sum _{k=m}^{n-r-2} \binom{n-r-k}{2}\log C_{k} (q) \}=0  \end{align*}
because of \eqref{PotF3}.
Finally, since
\begin{align*}   \{ \lan\; \ran ^{0+}_{1,0}  & +\frac{3(n-1-r-m)^2+n-r+m-3}{48} \log (1 - q\prod _{a=1}^r  l_a ^{l_a}\prod _{b=1}^m  (-l'_b) ^{l'_b}) \\
 & + \frac{1}{2}\sum _{k=m}^{n-r-2} \binom{n-r-k}{2}\log C_{k} (q) \} |_{q=0}=0,  \end{align*}
we are done.

\subsection{Corollaries}

First of all, if $m \geq 2$, then $I_0 =1$ and $ I_1^\zeta =0$.
Thus, we have

\begin{Cor}
If $m \geq 2$, then
\begin{align*}\lan \; \ran^\infty_{1,0} = \lan \; \ran^{0+}_{1,0} & =  -\frac{3(n-1-r-m)^2+n-r+m-3}{48} \log (1 - q\prod _{a=1}^r  l_a ^{l_a}\prod _{b=1}^m  (-l'_b) ^{l'_b}) \\
 & - \frac{1}{2}\sum _{k=m}^{n-r-2} \binom{n-r-k}{2}\log C_{k} (q) .
\end{align*}
\end{Cor}

If $m=1$, then $I_0 =1$ and 
\begin{align*}
\int_X H \cup c_{\text{dim}X -1}(T_X) = {n \choose 2} -\sum_{a=1}^r \frac{n}{l_a} - \frac{n}{l'_1}.
\end{align*}
Thus, we have

\begin{Cor}
If $m =1$, then
\begin{align*}\lan \; \ran^\infty_{1,0}|_{q \mapsto qe^{I_1(q)}}  &= \frac{I_1(q)}{24} \left( {n \choose 2} -\sum_a \frac{1}{l_a} -\frac{1}{l'_1} \right)  \\
&   -\frac{3(n-1-r-m)^2+n-r+m-3}{48} \log (1 - q\prod _{a=1}^r  l_a ^{l_a}\prod _{b=1}^m  (-l'_b) ^{l'_b}) \\
 & - \frac{1}{2}\sum _{k=m}^{n-r-2} \binom{n-r-k}{2}\log C_{k} (q) 
\end{align*}
where $I_1(q) \in \QQ[[q]]$ is the coefficient of $H$ in $I^\zeta_1$.
\end{Cor}

If $m=0$, then 
\begin{align*}
\int_X H \cup c_{\text{dim}X -1}(T_X) = {n \choose 2} -\sum_{a=1}^r \frac{n}{l_a} .
\end{align*}
Thus, we have

\begin{Cor}
If $m =0$, then
\begin{align*}\lan \; \ran^\infty_{1,0}|_{q \mapsto qe^{\frac{I_1(q)}{I_0(q)}}}  &= \frac{1}{24}\chi_{\mathrm{top}}(X)\log I_0 + \frac{1}{24}\frac{I_1(q)}{I_0(q)} \left( {n \choose 2} -\sum_a \frac{1}{l_a} \right)  \\
&   -\frac{3(n-1-r-m)^2+n-r+m-3}{48} \log (1 - q\prod _{a=1}^r  l_a ^{l_a}\prod _{b=1}^m  (-l'_b) ^{l'_b}) \\
 & - \frac{1}{2}\sum _{k=m}^{n-r-2} \binom{n-r-k}{2}\log C_{k} (q) 
\end{align*}
where $I_1(q) \in \QQ[[q]]$ is the coefficient of $H$ in $I^\zeta_1$.
\end{Cor}

\section{Example}

Let $X$ be the total space of $\mathcal{O}_{\PP^1}(-2) \boxtimes \mathcal{O}_{\PP^1}(-2)$ on $\PP^1 \times \PP^1$. Then it can be obtained by pull-back of $\mathcal{O}_{\PP^3}(-2)$ on $\PP^3$ under
Segre embedding $i : \PP^1 \times \PP^1 \hookrightarrow \PP^3.$
Note that the image of $i$ is a quadric hypersurface in $\PP^3$.
By using identification $H_2(X,\ZZ) \cong H_2(\PP^1 \times \PP^1, \ZZ) \cong \ZZ \times \ZZ$, define Gromov-Witten invariants
$$N_d := \sum_{d_1+d_2=d} \text{deg}[\overline{M}_{1,0}(X,(d_1,d_2))]^{\text{vir}}$$
for positive integer $d$.
In this section, we will calculate $N_d$, $d>0$, explicitly by using main theorem.

First of all, we apply Theorem 1.1 to $X$ by putting $n=3$, $r=1$, $m=1$ and $l_1=l'_1=2$ and obtain 
\begin{align}
\label{F2}
\langle \; \rangle^{0+}_{1,0}= -\frac{1}{12}\text{log}(1-16q) -\frac{1}{2}\text{log}(1+q\frac{q}{dq} \frac{I_1(q)}{I_0(q)})
\end{align}
where $I_0(q)$ and $I_1(q)$ are defined by
\begin{align}
\label{I2}
\sum_d q^d \frac{\prod_{k=1}^{2d}(2H+kz) \prod_{k=0}^{2d-1}(-2H-kz)}{\prod_{i=1}^4 \prod_{k=1}^d (H-\lambda_i +kz)} \equiv I_0(q) + I_1(q) \frac{H}{z} +O(\frac{1}{z^2}).
\end{align}
Here, left hand side of (\ref{I2}) is $I_\T(0,q)$ and '$\equiv$' means modulo (\ref{root}).
Precisely,
\begin{align}
\label{I01}
I_0(q) =1, ~~ \text{ and } ~~ I_1(q)= \sum_{d>0} \frac{q^d}{d}  {2d \choose d}^2.
\end{align}
By Theorem 1.3 and Corollary 3.6, we obtain
\begin{align}
\label{G2}
\langle \; \rangle^{\infty}_{1,0} |_{q^d \mapsto q^d\text{exp}(dI_1(q))}= \frac{1}{12} I_1(q)  +\langle \; \rangle^{0+}_{1,0}.
\end{align}
Let's define
$$Q:= q\text{exp}I_1(q)~~\text{ and }~~ T:=\text{log}Q.$$
Then, by combining (\ref{F2}), (\ref{I01}) and (\ref{G2}), we have
\begin{align}
\label{M1}
\langle \; \rangle^{\infty}_{1,0} (Q) &= \frac{T}{12} + \frac{1}{2} \text{log} \left( (1-16q)^{-\frac{1}{6}} q^{-\frac{7}{6}} \frac{dq}{dT}\right)\\
\label{M2}
&=-\frac{1}{3}q -\frac{11}{6}q^2 -\frac{124}{9}q^3 +O(q^4)
\end{align}
because
$$I_1(q)=T-\text{log}q=1+4q+18q^2+\frac{400}{3}q^3+O(q^4)$$ 
and 
$$1+q\frac{d}{dq}I_1(q)=q\frac{dT}{dq}=1+4q+36q^2+400q^3+O(q^4).$$

One remark we have to mention is that, in \cite{Hos}, he proves
$$-\frac{T}{12}+\langle \; \rangle^{\infty}_{1,0} (Q) = \frac{1}{2} \text{log} \left( (1-16q)^{-\frac{1}{6}} q^{-\frac{7}{6}} \frac{dq}{dT}\right)$$
which is exactly same as (\ref{M1}).

On the other hand, by using 
$$Q=q\text{exp}I_1(q)=q+4q^2+26q^3+O(q^4),$$ 
we have
\begin{align}
\label{M3}
\langle \; \rangle^{\infty}_{1,0}(Q) &= \sum_{d=1}^\infty  N_d Q^d \nonumber \\
&= N_1q+(4N_1 +N_2)q^2+(26N_1 + 8N_2 + N_3)q^3+O(q^4). 
\end{align}
Then, by comparing (\ref{M2}) and (\ref{M3}), we have
$$N_1=-\frac{1}{3},~~N_2=-\frac{1}{2} ~~ \text{ and }~~ N_3=-\frac{10}{9}, ~~\cdots.$$
Another remark is that these numbers are appeared in \cite{ABK}.
In this paper, they also showed that
\begin{align}
\label{modular} \langle \; \rangle^{\infty}_{1,0} (Q) = \frac{T}{12} -\log \eta(\tau)
\end{align}
where $\eta$ is Dedekind eta function and
$$Q=p^{\frac{1}{2}}-4p+6p^{\frac{3}{2}}+ \cdots ,~~p=e^{2\pi i \tau},$$
i.e,
$$\eta(\tau)=e^{\frac{\pi i \tau}{12}} \prod_{n=1}^\infty (1-e^{2n\pi i \tau}) =p^{\frac{1}{24}} \prod_{n=1}^\infty (1-p^n).$$
Indeed, mirror curves of $X$ is a family of elliptic curves. If we regard $\tau$ as a parameter of family of elliptic curves which corresponds to 
$$\CC / (\ZZ \oplus \tau \ZZ),$$
then they showed (\ref{modular}) by modular properties and behavior at the discriminant of the family of elliptic curves of partition function at genus 1 which is defined exactly same as $-\frac{T}{12}+\langle \; \rangle^{\infty}_{1,0} (Q).$


\end{document}